\documentclass[12pt,a4paper]{article}
\usepackage{amsfonts}
\usepackage{amsmath, amssymb, graphics, setspace}

\usepackage{url}

\newtheorem{obs}{Remark}[section]
\newtheorem{teor}{Theorem}[section]

\author{Ernest Scheiber\footnote{e-mail: scheiber@unitbv.ro}}
\title{A Convergence Theorem for the \textit{Parareal} Algorithm Revisited}
\date{}

\begin{document}
\maketitle

\begin{abstract}
The subject of the paper is to verify the convergence conditions for the parareal algorithm 
using Gander and Hairer's theorem . The analysis is conducted in the case where the coarse
integrator is the Euler method and the high-accuracy integrator is an explicit 
Runge-Kutta type method.
\\

2020 \textit{Mathematics Subject Classification:} 65M22, 65L05

\textit{Key words:} parareal algorithm,  initial value problem.
\end{abstract}

\section{Introduction}

\textit{Parareal} (parallel in real-time) is an iterative algorithm designed to solve the initial value problem (IVP):

\begin{eqnarray}
\dot{x}(t) &=& f(t,x(t)),\quad t\in [t_0,T],\label{prealet1}\\
x(t_0) &=& x_0, \label{prealet2}
\end{eqnarray}

where $f:[t_0,T]\times\mathbb{R}^d\rightarrow\mathbb{R}^d.$ 

Introduced in 2001 by Jacques-Louis Lions, Yvon Maday and Gabriel Turinici. \cite{2}, the parareal algorithm shares similarities with the multiple shooting method \cite{1}, \cite{6}. Moreover, it has found applications in solving partial differential equations.

The main appeal of the parareal algorithm lies in its capability for parallel execution. It combines two numerical methods for solving the IVP:

\begin{itemize}
\item A numerical method with reduced precision and computational operations;
\item A high-precision numerical method that demands a greater computational load.
\end{itemize}

The convergence of the algorithm has been extensively studied \cite{1}, \cite{3}, \cite{6}.

This paper aims to emphasize a convergence result whose conditions can be validated in concrete cases: the Euler method for the low-precision method and a four-level Runge-Kutta scheme for the high-precision method. A similar theme can be found in \cite{6}, \cite{4}, albeit in a different context.

The structure of the paper unfolds as follows: Section 2 elaborates on the parareal algorithm, while Section 3 revisits the convergence theorem alongside the verification of conditions for the Euler method and the four-level Runge-Kutta method pair.

\section{The Parareal Algorithm}

Let's assume that the problem (\ref{prealet1})-(\ref{prealet2}) possesses a unique solution.

Consider the mesh defined as:
\begin{equation}\label{prealet4}
t_0<t_1<\ldots <t_N=T
\end{equation}
and let $(u_n^k)_{0\le n\le N}$ denote the approximations of the solution of the 
IVP $(x(t_n))_{0\le n\le N}$ obtained at the $k$-th iteration. There will be $K$ iterations. Define $I_n=[t_{n-1},t_n],$ for $n\in{1,2,\ldots,N}$.

The Parareal algorithm consists of iteratively computing the layers $(u_n^{(k)})_{0\le n\le N}$, for $k=0,1,\ldots,K$.

We have two available methods, referred to as functions:

\begin{itemize}
\item $u_n^k=C_{I_n}(u_{n-1}^{(k)})$, a numerical integration method with a small number of operations and reduced precision (\textit{cheap coarse integrator});
\item $u_n^k=F_{I_n}(u_{n-1}^{(k)})$, a numerical integration method that offers a higher degree of precision (\textit{fine, high accuracy integrator}).
\end{itemize}

The index in $I_n$ specifies the interval in which the numerical integration of the system (\ref{prealet1}) takes place with the initial condition $x(t_{n-1})=u_{n-1}^{(k)}$. The numerical solution at $t_n$ for iteration $k$ is denoted by $u_n^{(k)}$. For instance,
\begin{itemize}
\item $u_n^{(k)}=C_{I_n}(u_{n-1}^{(k)})=u_{n-1}^{(k)}+(t_n-t_{n-1})f(t_{n-1},u_{n-1}^{(k)})$ is given by the Euler method;
\item $F_{I_n}(u_{n-1}^{(k)})$ represents one or several explicit Runge-Kutta steps with 4 stages ($m=4$).
\end{itemize}

The algorithm consists of two components:
\begin{enumerate}
\item Initialization:
\begin{equation}\label{prealet5}
u_0^{(0)}=x_0,\qquad u_n^{(0)}=C_{I_n}(u_{n-1}^{(0)}),\quad n=1,2,\ldots,N.
\end{equation}
\item Iterations: For $k=1,2,\ldots,K$, the following formulas are used:
\begin{eqnarray}
u_n^{(k)} &=& u_n^{(k-1)}, \quad n=0,1,\ldots, k-1; \label{prealet31} \\
u_n^{(k)} &= & C_{I_n}(u_{n-1}^{(k)})+F_{I_n}(u_{n-1}^{(k-1)})-C_{I_n}(u_{n-1}^{(k-1)}), \quad n=k,k+1,\ldots,N. \nonumber\\
&& \hspace{6cm}\label{prealet3}
\end{eqnarray}
\end{enumerate}

$$
\begin{array}{c|c|c|c|c|c|c|}
& I_1 & I_1 & I_3 & \ldots & I_{N-1} & I_N \\
\hline
\mbox{Initialization}& u_0^{(0)}=x_0 & u_1^{(0)} & u_2^{(0)} & \ldots & u_{N-1}^{(0)} & u_N^{(0)} \\
\hline
k=1 &  \downarrow & u_1^{(1)} & u_2^{(1)} & \ldots & u_{N-1}^{(1)} & u_N^{(1)} \\
\hline
k=2 &  \downarrow & \downarrow & u_2^{(2)} & \ldots & u_{N-1}^{(2)} & u_N^{(2)} \\
\hline
 &  \downarrow & \downarrow & \downarrow &  &  & \vdots \\
\hline
\end{array}
$$

\begin{obs}\end{obs}
For $n=k,k+1,\ldots,N$, the quantity $\xi_n=F_{I_n}(u_{n-1}^{(k-1)})-C_{I_n}(u_{n-1}^{(k-1)})$ can be computed in parallel.

\begin{obs}\end{obs}
Let the sequence $u_n=F_{I_n}(u_{n-1}), n\in{1,2,\ldots,N}, u_0=x_0.$ Then:
\begin{itemize}
\item $u^{(k)}_0=u_0=x_0$ for all $k\in\{1,2,\ldots,K\}$.
\item The following equalities hold: $u^{(k)}_n=u_n$ for any $n\in\{1,2,\ldots,k\}$.

Proof by induction on $k$:

1. $k=1$

\noindent
The following equalities hold:
$$
u^{(1)}_1=C_{I_1}(u^{(1)}_0)+F_{I_1}(u^{(0)}_0)-C_{I_1}(u^{(0)}_0)=C_{I_1}(x_0)+F_{I_1}(x_0)-C_{I_1}(x_0)=
$$
$$
=F_{I_1}(x_0)=u_1.
$$
From (\ref{prealet31}) it follows that $u^{(k)}_1=u_1, \forall\ k\ge 1.$ 

2. $u^{(k-1)}_{k-1}=u^{(k)}_{k-1}=u_{k-1} \quad\Rightarrow\quad u^{(k)}_{k}=u_{k}.$

Indeed
$$
u^{(k)}_{k}=C_{I_{k}}(u^{(k)}_{k-1})+F_{I_{k}}(u^{(k-1)}_{k-1})-C_{I_{k}}(u^{(k-1)}_{k-1})=
$$
$$
=C_{I_{k}}(u_{k-1})+F_{I_{k}}(u_{k-1})-C_{I_{k}}(u_{k-1})=F_{I_{k}}(u_{k-1})=u_{k}.
$$

The recursive formula (\ref{prealet3}) also holds for $n\in\{1,2,\ldots,k-1\}$.
\end{itemize}

\begin{obs}\end{obs}
When $K=N$, the solution provided by the parareal algorithm is the same as that given by the high accuracy integrator. 
Practically, the only interesting case is when $K<N$.

\section{Gander and Hairer's convergence theorem}

We shall follow the presentation of the convergence theorem given by M. J. Gander \c{s}i E. Hairer \cite{1}.

\noindent
Let be 
\begin{itemize}
\item
$\|\cdot\|$ be a norm in $\mathbb{R}^d;$ 
\item
$ h=\frac{T-t_0}{N}$ which implies $t_n=t_0+n h,\ n\in\{0,1,\ldots,N\};$
\end{itemize}

\begin{teor}\label{prealt1}
Let $(u_n)_{0\le n\le N}$ be the numerical solution of the problem (\ref{prealet1})-(\ref{prealet2})
given by the high accuracy integrator, $u_n=F_{I_n}(u_{n-1}),
n\in\{1,2,\ldots,N\}.$ If
\begin{enumerate}
\item 
$$
\|C_{I_n}(u_1)-C_{I_n}(u_2)\|\le \underbrace{(1+h c_1)}_{b}\|u_1-u_2\|,
\qquad \forall\ u_1,u_2\in \mathbb{R}^d, c_1>0;
$$
\item
$$
\|F_{I_n}(u)-C_{I_n}(u)\|\le h^{1+\alpha}c_2,\qquad \forall u\in \mathbb{R}^d,\ \alpha>0,\ c_2>0;
$$
\item
$$
\|(F_{I_n}(u_1)-C_{I_n}(u_1))-(F_{I_n}(u_2)-C_{I_n}(u_2))\|\le \underbrace{h c_3}_{a}\|u_1-u_2\|,
$$
$
\hspace*{4cm} \forall\ u_1,u_2\in \mathbb{R}^d, c_3>0
$
\end{enumerate}
then
$$
\lim_{h\searrow 0}\|u^{(k)}-u\|_h=\lim_{h\searrow 0}\ \max_{0\le n\le N}\|u_n^{(k)}-u_n\|=0.
$$
\end{teor}

\vspace*{0.3cm}\noindent
\textbf{Proof.}
We shall use the notations
$$
E_n^{(k)}=\|u_n^{(k)}-u_n\|,\quad n\in\{0,1,\ldots,N\}.
$$

For $n\in\{1,2,\ldots,N\}$ \c{s}i $k\ge 1,$ from the recurrence (\ref{prealet3}) we obtain
$$
u_n^{(k)}-u_n=
$$
$$
=(C_{I_n}(u_{n-1}^{(k)})-C_{I_n}(u_{n-1}))+((F_{I_n}(u_{n-1}^{(k-1)})-C_{I_n}(u_{n-1}^{(k-1)}))-
(F_{I_n}(u_{n-1})-C_{I_n}(u_{n-1})))
$$
The hypotheses of the theorem imply
$$
\|u_n^{(k)}-u_n\|\le
$$
$$
\le\|C_{I_n}(u_{n-1}^{(k)})-C_{I_n}(u_{n-1}))\|+\|(F_{I_n}(u_{n-1}^{(k-1)})-C_{I_n}(u_{n-1}^{(k-1)}))-
(F_{I_n}(u_{n-1})-C_{I_n}(u_{n-1}))\|\le
$$
$$
\le b\|u_{n-1}^{(k)}-u_{n-1}\|+a\|u_{n-1}^{(k-1)}-u_{n-1}\|.
$$
The above inequality may be rewritten as
\begin{equation}\label{prealet7}
E_n^{(k)}\le b E_{n-1}^{(k)}+a E_{n-1}^{(k-1)}.
\end{equation}

\noindent
For $k=0$ we obtain
$$
u^{(0)}_n-u_n=C_{I_n}(u^{(0)}_{n-1})-F_{I_n}(u_{n-1})=
$$
$$
=(C_{I_n}(u^{(0)}_{n-1})-C_{I_n}(u_{n-1}))+ (C_{I_n}(u_{n-1})-  F_{I_n}(u_{n-1})).
$$
It results that
$$
\|u^{(0)}_n-u_n\|\le \|C_{I_n}(u^{(0)}_{n-1})-C_{I_n}(u_{n-1})\|+\|C_{I_n}(u_{n-1})-  F_{I_n}(u_{n-1})\|\le
$$
$$
\le b\|u^{(0)}_{n-1}-u_{n-1}\|+\underbrace{h^{1+\alpha}c_2}_{\gamma},
$$
thus
\begin{equation}\label{prealet8}
E^{(0)}_{n}\le b E^{(0)}_{n-1}+\gamma.
\end{equation}

Its give the idea to study the sequence  $(z^{(k)}_n)_{n,k\in\mathbb{N}}$ defined by the recurrence formulas
$$
\begin{array}{lclcl}
z^{(k)}_0 &=& 0, & \qquad& k\in \{0,1,\ldots\};\\
z^{(0)}_n &=& b z^{(0)}_{n-1}+\gamma, &\qquad & n\in\{1,2,\ldots\};\\
z_n^{(k)} &=& b z_{n-1}^{(k)}+a z_{n-1}^{(k-1)}, &\qquad &  n\in\{1,2,\ldots\},\ k\ge 1.
\end{array}
$$
We retain the inequality $E^{(k)}_n\le z^{(k)}_n.$

The \textit{generating} function  $\rho_k(\zeta)=\sum_{n\ge 1}z^{(k)}_n\zeta^n$ verifies the equalities
$$
\begin{array}{lclcl}
\rho_k(\zeta) &=& a\zeta\rho_{k-1}(\zeta)+b\zeta\rho_k(\zeta) & \Rightarrow & \rho_k(\zeta)=\frac{a\zeta}{1-b\zeta}\rho_{k-1}(\zeta);\\
\rho_0(\zeta) &=& b\zeta\rho_0(\zeta)+\frac{\gamma\zeta}{1-\zeta} & \Rightarrow & \rho_0(\zeta)=\frac{\gamma\zeta}{(1-\zeta)(1-b\zeta)}.
\end{array}
$$
We find
$$
\rho_k(\zeta)=\left(\frac{a\zeta}{1-b\zeta}\right)^k\rho_0(\zeta)=\frac{\gamma a^k}{1-\zeta}\left(\frac{\zeta}{1-b\zeta}\right)^{k+1}.
$$
Because $b>1,$ for $0<\zeta<\frac{1}{b}$ it results $\frac{1}{1-\zeta}\le \frac{1}{1-b\zeta},$ and then 
$$
\rho_k(\zeta)\le \frac{\gamma a^k \zeta^{k+1}}{(1-b\zeta)^{k+2}}.
$$
Using serial expansion
$$
\frac{1}{(1-b\zeta)^{k+1}}=\sum_{m=0}^{\infty}\left(
\begin{array}{c}m+k\\k\end{array}\right)b^m\zeta^m.
$$
it will result 
$$
\frac{\gamma a^k \zeta^{k+1}}{(1-b\zeta)^{k+2}}=\gamma a^k \sum_{m=0}^{\infty}\left(
\begin{array}{c}m+k+1\\k+1\end{array}\right)b^m\zeta^{m+k+1}.
$$
The coefficient of  $\zeta^n$ is obtained for $m=n-k-1$ and it is
$$
 \gamma a^k \beta^{n-k-1}\left(
\begin{array}{c}n\\k+1\end{array}\right).
$$
It results the inequality
$$
E^{(k)}_n\le z^{(k)}_n\le \gamma a^k \beta^{n-k-1}\frac{n(n-1)\ldots(n-k)}{(k+1)!}=
$$
$$
=h^{1+\alpha}c_2 (h c_3)^k (1+h c_1)^{n-k-1}\frac{n(n-1)\ldots(n-k)}{(k+1)!}\le
$$
$$
\le h^{\alpha}\frac{c_2c_3^k (T-t_0)^{k+1} e^{c_1 (T-t_0)}}{(k+1)!}. \quad\rule{5pt}{5pt}
$$

\subsection*{An application}
We shall verify the hypotheses of the above theorem when Euler method is the coarse integrator $C_{I_n}$ and 
when the high accuracy integrator is the Runge-Kutta method with four levels $F_{I_n}.$
In this case
\begin{eqnarray*}
C_{I_n}(u)&=&u+hf(t_{n-1},u) \\
F_{I_n}(u)&=& u+hF_4(h,t_{n-1},u;f)\\
&& \mbox{where}\\
&& F_4(h,t,u;f)=\frac{1}{6}(k_1(h)+2k_2(h)+2k_3(h)+k_4(h));\\
&& k_1(h)=f(t,u) \\
&& k_2(h)=f(t+\frac{h}{2},u+\frac{h}{2}k_1(h))\\
&& k_3(h)=f(t+\frac{h}{2},u+\frac{h}{2}k_2(h))\\
&& k_4(h)=f(t+h,u+hk_3(h))
\end{eqnarray*}

We assume that the function $f$ satisfies the Lipschitz condition
$$ \|f(t,u_1)-f(t,u_2)\|\le L\|u_1-u_2\|, \ \forall\ u_1,u_2\in \mathbb{R}^d,\ \forall\ t\in[t_0,T].$$
This assumption implies the existence and the bounding of the IVP as well as the 
bounding of the numerical solution of any convergent numerical method.

\begin{enumerate}
\item
The equality
$$
C_{I_n}(u_1)-C_{I_n}(u_2)=u_1-u_2+h(f(t_{n-1},u_1)-f(t_{n-1},u_2))
$$
implies
$$
\|C_{I_n}(u_1)-C_{I_n}(u_2)\|\le (1+hL)\|u_1-u_2\|.
$$ 
That is the first condition with $c_1=L.$
\item
The following equality holds
\begin{equation}\label{obset1}
F_{I_n}(u)-C_{I_n}(u)=h\left(F_4(h,t_{n-1},u;f)-f(t_{n-1},u)\right).
\end{equation}
If the function $f$ is smooth enough then the\textit{Mathematica} code
\begin{doublespace}
\noindent\(\pmb{\text{K1}=f[t,u];}\)
\end{doublespace}

\begin{doublespace}
\noindent\(\pmb{\text{K2}=\text{Series}[f[t+h/2,u+h/2 \text{K1}],\{h,0,1\}];}\)
\end{doublespace}

\begin{doublespace}
\noindent\(\pmb{\text{K3}=\text{Series}[f[t+h/2,u+h/2 \text{K2}],\{h,0,1\}];}\)
\end{doublespace}

\begin{doublespace}
\noindent\(\pmb{\text{K4}=\text{Series}[f[t+h,u+h \text{K3}],\{h,0,1\}];}\)
\end{doublespace}

\begin{doublespace}
\noindent\(\pmb{1/6(\text{K1}+2 \text{K2}+2 \text{K3}+\text{K4})-f[t,u]}\)
\end{doublespace}

\begin{doublespace}
\noindent\(\frac{1}{6} \left(3 f[t,u] f^{(0,1)}[t,u]+3 f^{(1,0)}[t,u]\right) h+O[h]^2\)
\end{doublespace}

\begin{doublespace}
\noindent\(\pmb{\text{Simplify}[\text{Collect}[1/6(\text{K1}+2 \text{K2}+2 \text{K3}+\text{K4})-f[t,u],h]]}\)
\end{doublespace}

\begin{doublespace}
\noindent\(\frac{1}{2} h \left(f[t,u] f^{(0,1)}[t,u]+f^{(1,0)}[t,u]\right)\)
\end{doublespace}

\noindent
proves that the expression in the brackets from (\ref{obset1}) is of the form 
$\Phi_1(t,u) h+o(h^2).$ Consequently
\begin{equation}\label{prealet10}
\|F_{I_n}(u)-C_{I_n}(u)\|\ \le h^2 \Lambda,
\end{equation}
where $\Lambda$ is an upper bound of  $\Phi_1(t,u)$ in a compact set that includes the graph of the solution in the interval  $[t_0,T].$
\item
The function $F_4(h,t,x;f)$ satisfies the Lipschitz condition, too, 
$$,
\|F_4(h,t,u;f)-F_4(h,t,v;f)\|\le M\|u-v\|,
$$
with $M=L\left(1+\frac{1}{2}hL+\frac{1}{6}h^2L^2+\frac{1}{24}h^3L^3\right).$
It follows that
$$
\|(F_{I_n}(u)-C_{I_n}(u))-(F_{I_n}(v)-C_{I_n}(v))\|\le h(L+M)\|u-v\|.
$$
\end{enumerate}

Thus, we find
\begin{enumerate}
\item 
$b:=1+hL,\ c_1:=L;$
\item
$\alpha:=1;$
\item
$a:=h c_3,\ c_3:=L+M.$
\end{enumerate}

In the above version, the action of the method $F_{I_n}$ consists of a single Runge-Kutta step on the interval of length $h.$
We are concerned with the variant in which a number of $m$ Runge-Kutta steps are executed on intervals of length $\tau=h/m.$

\begin{enumerate}
\setcounter{enumi}{1}
\item
Let us denote $t_{n,j}=t_{n-1}+j\tau,\ j\in\{0,1,\ldots,m\}$ and $I_{n,j}=[t_{n,j-1},t_{n,j}].$
Now, we shall suppose that the differential system  (\ref{prealet1}) is autonomous,
$
\dot{x}(t)=f(x(t)),
$
satisfying the Lipschitz condition.

We begin computing the function $F_{I_n}(u).$
Writing $\hat{u}_0=u$ the following equalities occur
\begin{eqnarray*}
\hat{u}_1 &=&\hat{u}_0+\tau F_4(\tau,t_{n,0},\hat{u}_0;f)\\
\hat{u}_2 &=& \hat{u}_1+\tau F_4(\tau,t_{n,1},\hat{u}_1;f)\\
\vdots\\
\hat{u}_m &=& \hat{u}_{m-1}+\tau F_4(\tau,t_{n,m-1},\hat{u}_{m-1};f)
\end{eqnarray*}
and consequently
$$
F_{I_n}(u)=\hat{u}_m=u+\tau \sum_{j=0}^{m-1}F_4(\tau,t_{n,j},\hat{u}_j;f).
$$

\noindent
Then we have
\begin{equation}\label{prealet9}
F_{I_n}(u)-C_{I_n}(u)=\tau \sum_{j=0}^{m-1}F_4(\tau,t_{n,j},\hat{u}_j;f)-hf(t_{n-1},u)=
\end{equation}
$$
=\tau \sum_{j=0}^{m-1}\left(F_4(\tau,t_{n,j},\hat{u}_j;f)-f(t_{n-1},u)\right)=
$$
\begin{equation}\label{prealet11}
=\tau \sum_{j=0}^{m-1}\left(F_4(\tau,t_{n,j},\hat{u}_j;f)-f(t_{n,j},\hat{u}_j)\right)+
\tau \sum_{j=0}^{m-1}(f(t_{n,j},\hat{u}_j)-f(t_{n,0},\hat{u}_0)).
\end{equation}

Taking into account the justification of inequality (\ref{prealet10}), it follows that
\begin{equation}\label{prealet15}
\| \tau ( F_4(\tau,t_{n,j},\hat{u}_j;f)-f(t_{n,j},\hat{u}_j))\|\le \tau^2\Lambda.
\end{equation}

We will proceed to establish an upper bound for $\|\hat{u}_j-u\|,\ j\in\{1,2,\ldots,m\}.$
The inequality occurs
\begin{equation}\label{prealet12}
\|\hat{u}_j-u\| \le \|\hat{u}_j-\hat{u}(t_{n,j})\| + \|\hat{u}(t_{n,j})-u\|.
\end{equation}

\begin{enumerate}
\item
Let $\hat{u}_{\tau}=(\hat{u}_j)_{0\le j\le m}$ be the  solution of the IVP
\begin{eqnarray*}
\dot{x}(t) &=& f(x) \\
x(t_{n,0}) &=& u
\end{eqnarray*}
and let $[\hat{u}]_{\tau}=(\hat{u}(t_{n,0}),\hat{u}(t_{n,1}),\ldots,\hat{u}(t_{n,m}))$ 
be the numerical solution computed using the Runge-Kutta method with four levels.
Based on the consistency and stability \cite{10}, the following inequality occurs
$$
 \|\hat{u}_{\tau}-[\hat{u}]_{\tau}\|_{\tau}\le c_4\tau^4,
$$
where $\|\cdot\|_{\tau}$ is the maximum norm in $\{(u_k)_{0\le k\le m} : u_k\in\mathbb{R}^d\}.$
Thus $\|\hat{u}_j-\hat{u}(t_{n,j})\|\le c_4\tau^4,\ \forall\ j\in\{0,1,\ldots,m\}.$
\item
From
$$
\hat{u}(t_{n,j})-u=\hat{u}(t_{n,j})-\hat{u}(t_{n,0})=\int_{t_{n,0}}^{t_{n,j}}f(\hat{u}(s))\mathrm{d}s=
$$
$$
=\int_{t_{n,0}}^{t_{n,j}}(f(\hat{u}(s))-f(\hat{u}(t_{n,0}))\mathrm{d}s+f(\hat{u}(t_{n,0}))j\tau
$$
we deduce
$$
\|\hat{u}(t_{n,j})-u\|\le L\int_{t_{n,0}}^{t_{n,j}}\|\hat{u}(s))-u)\|\mathrm{d}s+j\tau \|f(u)\|.
$$
Using the Gr\"{o}nwall's Lemma it results
$$
\|\hat{u}(t_{n,j})-u\|\le \underbrace{\|f(u)\|e^{Lh}}_{\le\  c_5}j\tau \le c_5 h.
$$
\end{enumerate}
Based on  (\ref{prealet12}) it results the upper bound
\begin{equation}\label{prealet14}
\|\hat{u}_j-u\| \le c_4\tau^4+ c_5 h.
\end{equation}
and from (\ref{prealet11}) we obtain
$$
\|F_{I_n}(u)-C_{I_n}(u)\|\le \tau^2\Lambda m +m \tau L ( c_4\tau^4+ c_5 h)=
h^2\underbrace{\left(\frac{\Lambda}{m}+\frac{L c_4 h^3}{m^4}+L c_5\right)}_{\tilde{\Lambda}}.
$$
\item
For any  $u,v\in \mathbb{R}^d$ 
\begin{eqnarray*}
F_{I_n}(u) &=& u+\tau \sum_{j=0}^{m-1}F_4(\tau,t_{n,j},\hat{u}_j;f);\\
F_{I_n}(v) &=& v+\tau \sum_{j=0}^{m-1}F_4(\tau,t_{n,j},\hat{v}_j;f).
\end{eqnarray*}
and then
$$
F_{I_n}(u)-C_{I_n}(u))-(F_{I_n}(v)-C_{I_n}(v))=
$$
\begin{equation}\label{prealet13}
=\tau \sum_{j=0}^{m-1}(F_4(\tau,t_{n,j},\hat{u}_j;f)-F_4(\tau,t_{n,j},\hat{v}_j;f))-h(f(u)-f(v)).
\end{equation}
The equality
$$
\hat{u}_j -\hat{v}_j=\hat{u}_{j-1} -\hat{v}_{j-1}+\tau (F_4(\tau,t_{n,-1},\hat{u}_{j-1};f)-F_4(\tau,t_{n,-1},\hat{v}_{j-1};f))
$$
and the Lipschitz condition of $F_4(h,t,x;f)$ implies
$$
\|\hat{u}_j -\hat{v}_j\| \le (1+\tau M)\|\hat{u}_{j-1} -\hat{v}_{j-1}\|,\qquad j\in\{1,2,\ldots,m-1\}.
$$
It results the inequality
$$
\|\hat{u}_j -\hat{v}_j\| \le (1+\tau M)^j\|u-v\|
$$
and then
$$
\left\| \tau \sum_{j=0}^{m-1}(F_4(\tau,t_{n,j},\hat{u}_j;f)-F_4(\tau,t_{n,j},\hat{v}_j;f)) \right\| \le
\tau M \sum_{j=0}^{m-1} (1+\tau M)^j \|u-v\| =
$$
$$
= ((1+\tau M)^m-1) \|u-v\|\ \le (e^{hM}-1) \|u-v\| = h \underbrace{\frac{e^{hM}-1}{h}}_{\le\ c_6} \|u-v\|\le hc_6 \|u-v\|,
$$
because $\lim_{h\rightarrow 0}\frac{e^{hM}-1}{h}=M.$
From(\ref{prealet13}) we find
$$
\|F_{I_n}(u)-C_{I_n}(u))-(F_{I_n}(v)-C_{I_n}(v))\|\le h\underbrace{(c_6+L)}_{\tilde{c}_3}\|u-v\|.
$$
\end{enumerate}

Evaluations similar to those deduced for $m=1$ have been derived, with which
the convergence conditions from Theorem \ref{prealt1} were verified.

\vspace*{0.3cm}
Hypothesis 3 of the Theorem \ref{prealt1} can be dropped:
\begin{teor}\label{prealt2}
Let $(u_n)_{0\le n\le N}$ be the numerical solution of the problem (\ref{prealet1})-(\ref{prealet2})
given by the high accuracy integrator, $u_n=F_{I_n}(u_{n-1}),
n\in\{1,2,\ldots,N\}.$ If
\begin{enumerate}
\item 
$$
\|C_{I_n}(u_1)-C_{I_n}(u_2)\|\le \underbrace{(1+h c_1)}_{b}\|u_1-u_2\|,
\qquad \forall\ u_1,u_2\in \mathbb{R}^d, c_1>0;
$$
\item
$$
\|F_{I_n}(u)-C_{I_n}(u)\|\le h^{1+\alpha}c_2,\qquad \forall u\in \mathbb{R}^d,\ \alpha>0,\ c_2>0;
$$
\end{enumerate}
then
$$
\lim_{h\searrow 0}\|u^{(k)}-u\|_h=\lim_{h\searrow 0}\ \max_{0\le n\le N}\|u_n^{(k)}-u_n\|=0.
$$
\end{teor}

\vspace*{0.3cm}\noindent
\textbf{Proof.}
With the above introduced notations, from the equality
$$
u_n^{(k)}-u_n=
$$
$$
=(C_{I_n}(u_{n-1}^{(k)})-C_{I_n}(u_{n-1}))+((F_{I_n}(u_{n-1}^{(k-1)})-C_{I_n}(u_{n-1}^{(k-1)}))-
(F_{I_n}(u_{n-1})-C_{I_n}(u_{n-1})))
$$
we deduce
$$
\|u_n^{(k)}-u_n\|\le
$$
$$
\le\|C_{I_n}(u_{n-1}^{(k)})-C_{I_n}(u_{n-1}))\|+
$$
$$
+\|(F_{I_n}(u_{n-1}^{(k-1)})-C_{I_n}(u_{n-1}^{(k-1)}))\|+
\|(F_{I_n}(u_{n-1})-C_{I_n}(u_{n-1}))\|\le
$$
$$
\le b\|u^{(k)}_{n-1}-u_{n-1}\|+2h^{1+\alpha}c_2.
$$
The last relation may be rewritten as
$$
E^{(k)}_n\le b E^{(k)}_{n-1}+2h^{1+\alpha}c_2.
$$
It results that
$$
E^{(k)}_n\le 2 h^{1+\alpha}c_2(1+b+\ldots+b^{n-1})\le 2h^{\alpha}\frac{e^{(T-t_0)c_1}c_2}{c_1}.\quad\rule{5pt}{5pt}
$$

With this version we may verify the parareal algorithm convergence when $C_{I_n}$
uses the backward Euler method.
$$
C_{I_n}(u)=z\qquad\mbox{where }z\mbox{ is the solution of the equation }z=u+h f(t_n,z).
$$
We verify the conditions of Theorem ref{prealt2} when the system is autonomous and the function $f$ satisfies the Lipschitz
condition.
\begin{enumerate}
\item
For $C_{I_n}(u_i)=z_i$ with $z_i=u_i+hf(z_i),\ i=1,2,$ the equalities occur
$$
C_{I_n}(u_1)-C_{I_n}(u_2)=z_1-z_2=u_1-u_2-h(f(z_1)-f(z_2)).
$$
It deduces that
$$
\|z_1-z_2\|\le\|u_1-u_2\|+hL\|z_1-z_2\|\quad\Leftrightarrow\quad \|z_1-z_2\|\le \frac{1}{1-hL}\|u_1-u_2\|
$$
Thus
$$
\|C_{I_n}(u_1)-C_{I_n}(u_2)\|\le \frac{1}{1-h L}\|u_1-u_2\|=\left(1+\frac{h L}{1-h L}\right)\|u_1-u_2\|.
$$
If $h\le \frac{1}{2L},$ in addition we have  $\frac{1}{1-h L}\le 2,$ and then
$$
\|C_{I_n}(u_1)-C_{I_n}(u_2)\|\le (1+2hL)\|u_1-u_2\|.
$$
\item
With the used notations, the following equalities occur
$$
F_{I_n}(u)-C_{I_n}(u)=u+\tau \sum_{j=0}^{m-1}F_4(\tau,t_{n,j},\hat{u}_j;f)-z=
$$
$$
=\tau  \sum_{j=0}^{m-1}(F_4(\tau,t_{n,j},\hat{u}_j;f)-f(\hat{u}_j))+\tau \sum_{j=0}^{m-1}(f(\hat{u}_j)-f(z)).
$$
Using  (\ref{prealet15}) we have
$$
\| \tau ( F_4(\tau,t_{n,j},\hat{u}_j;f)-f(t_{n,j},\hat{u}_j))\|\le \tau^2\Lambda.
$$
and from (\ref{prealet14}) we find
$$
\|\hat{u}_j-z\|\le \|\hat{u}_j-u\|+\|u-z\|\le c_4\tau^4+c_5h+hL\|f(z)\|.
$$
The numerical solution is  bounded and the Lipschitz condition of $f$ implies that it is is bounded on the set of numerical solution.
Consequently
$$
\|F_{I_n}(u)-C_{I_n}(u)\|\le m\tau^2\Lambda+m \tau L(c_4\tau^4+c_5h+hL\|f(z)\|)\le
$$
$$
\le h^2\left(\underbrace{\frac{\Lambda}{m}+\frac{Lc_4\tau^3}{m^3}+Lc_5+L^2\|f(z)\|}_{\le\tilde{c}_2}\right).
$$
\end{enumerate}
Thus
\begin{enumerate}
\item
$b:=1+2hL,\quad c_1:=2L;$
\item
$\alpha:=1,\quad c_2:=\tilde{c}_2.$
\end{enumerate}

\begin{obs}\end{obs} 
If the numerical methods defined by $C_{I_n}$ (the coarse integrator) and $F_{I_n}$ (the fine integrator)
have a convergence order greater than or equal to 2, then the second condition of Theorem \ref{prealt2} is satisfied.

\input cyracc.def
\font\tencyr=wncyr10
\def\cyr{\tencyr\cyracc}

\end{document}